
\nopagenumbers               
\footnote{{}}{\hskip -2 mm 2000 {\it Mathematics Subject
Classification}. Primary 11M06, Secondary 11F72.}

\def\txt#1{{\textstyle{#1}}}
\baselineskip=13pt
\def\hf{{\textstyle{1\over2}}}
\def\a{\alpha}
\def\d{{\,\rm d}}
\def\e{\varepsilon}
\def\f{\varphi}
\def\G{\Gamma}

\def\k{\kappa}
\def\s{\sigma}

\def\={\;=\;}
\def\zx{\zeta(\hf+ix)}
\def\zt{\zeta(\hf+it)}

\def\o{\omega}  
\def\R{\Re{\rm e}\,}  \def\s{\sigma}
\def\z{\zeta}

\def\no{\noindent} \def\o{\omega}
\def\H{H_j^3({\txt{1\over2}})}  \def\={\,=\,}
\def\hf{{\textstyle{1\over2}}}
\def\txt#1{{\textstyle{#1}}}
\def\f{\varphi}
\def\Z{{\cal Z}}
\font\tenmsb=msbm10
\font\sevenmsb=msbm7
\font\fivemsb=msbm5
\newfam\msbfam
\textfont\msbfam=\tenmsb
\scriptfont\msbfam=\sevenmsb
\scriptscriptfont\msbfam=\fivemsb
\def\Bbb#1{{\fam\msbfam #1}}

\def \NN {\Bbb N}
\def \CC {\Bbb C}
\def \RR {\Bbb R}
\def \ZZ {\Bbb Z}

\font\aa=cmcsc10 at 12pt 
\font\bb=cmcsc10
\font\cc=cmcsc10 at 8pt 

\font\ff=cmr9
\font\hh=cmbx12 
\def\rightheadline{\hfil\ff  Estimation of $\Z_2(s)$\hfil\folio}
\def\leftheadline{\ff\folio\hfil Aleksandar Ivi\'c
 \hfil}
\def\emptyheadline{\hfil}
\headline{\ifnum\pageno=1 \emptyheadline\else
\ifodd\pageno \rightheadline \else \leftheadline\fi\fi}
\topglue2cm
\centerline{\aa  ON THE ESTIMATION OF $\Z_2(s)$}
\bigskip\bigskip
\centerline{\aa Aleksandar Ivi\'c}
\bigskip
\centerline{\sevenbf Anal. Probab. Number Theory (A. Dubickas et al. eds.),
TEV, 2002, Vilnius, 83-98}
\bigskip
\centerline{\bb Abstract}
\bigskip\bigskip\no
{\ff Estimates for $\Z_2(s) = \int_1^\infty|\zx|^4x^{-s}\d x\;(\R s > 1)$  
are discussed, both pointwise and in the mean square. It is shown how these
estimates can be used to bound $E_2(T)$, the error term in the asymptotic
formula for $\int_0^T|\zt|^4\d t$.}

\bigskip\bigskip
{\hh 1. Introduction}
\bigskip

\noindent

The function $\Z_2(s)$  is the analytic continuation of the function
$$
\Z_2(s) \= \int_1^\infty |\zx|^4x^{-s}\d x \qquad(\R s > 1),
$$
and represents the (modified) Mellin transform of $|\zx|^4$.
It was introduced by Y. Motohashi
[15] (see also  [7], [10], [11]  and [16]), who showed  that  it  
has meromorphic continuation over $\CC$. 
In the half-plane $\s = \Re{\rm e}\, s > 0$
it has the following singularities: the pole $s = 1$ of order five, 
simple poles at $s = {1\over2} \pm i\k_j\,(\k_j
= \sqrt{\lambda_j - {1\over4}})$ and poles at $s = \hf\rho$, 
where $\rho$ denotes complex
zeros of $\zeta(s)$.  Here as usual $\,\{\lambda_j = \k_j^2 +
{1\over4}\} \,\cup\, \{0\}\,$ is the discrete spectrum of the
non-Euclidean Laplacian acting on $SL(2,\ZZ)$-automorphic forms
(see [16, Chapters 1--3] for a 
comprehensive account of  spectral theory and the Hecke $L$-functions).

\medskip
The  aim of this note is to study the estimation $\Z_2(s)$, both 
pointwise and in mean square. This  research was begun in [11], and 
continued in [7]. It was shown there that we have
$$
\int_0^T|{\cal Z}_2(\s + it)|^2\d t \ll_\e 
T^\e\left(T + T^{2-2\s\over1-c}\right) \qquad(\hf < \s < 1),\leqno(1.1)
$$
and we also have unconditionally
$$
\int_0^T|\Z_2(\s + it)|^2\d t \;\ll\;T^{10-8\s\over3}\log^CT
\qquad(\hf < \s <1,\,C > 0).\leqno(1.2)
$$
Here and later $\e$ denotes arbitrarily small, positive constants,
which are not necessarily the same ones at each occurrence, while
$\s$ is assumed to be fixed. The constant $c$ 
appearing in (1.1) is defined by
$$
E_2(T) \ll_\e T^{c+\e},\leqno(1.3)
$$
where the function $E_2(T)$ denotes the error term in the asymptotic
formula for the mean fourth power of $|\zt|$. It is customarily
defined by the relation
$$
\int_0^T|\zt|^4\d t \;=\; TP_4(\log T) \;+\;E_2(T),\leqno(1.4)
$$
with
$$
P_4(x) \;=\; \sum_{j=0}^4\,a_jx^j, \quad a_4 = {1\over2\pi^2}.\leqno(1.5)
$$
For the explicit evaluation of the $a_j$'s in (1.5), see [3].

Mean value estimates for $\Z_2(s)$ are a natural tool to investigate
the eighth power moment of $|\zt|$. Indeed, one has (see [7, (4.7)])
$$
\int_T^{2T}|\zt|^8\d t \ll_\e 
T^{2\s-1}\int_0^{T^{1+\e}}|\Z_2(\s+it)|^2\d t \quad(\hf < \s < 1).
\leqno(1.6)
$$ 
We shall prove here the pointwise estimate for $\Z_2(s)$ given by

\bigskip
THEOREM 1. {\it For $\hf < \s \le 1$ fixed and $t \ge t_0 > 0$ we have}
$$
\Z_2(\s + it) \;\ll_\e\; t^{{4\over3}(1-\s)+\e}.\leqno(1.7)
$$

\bigskip\no
Theorem 1 corrects an oversight in the proof of Theorem 3 of [7], where
the better exponent $1-\s$ was claimed, since unfortunately the
condition $1/3 \le \xi \le 1/2$ (see [11, (4.22)]) has to be observed,
and our argument needs $\xi = \e$ to hold. 
It improves the exponent $2(1-\s)$  that was obtained in [11].
Probably the exponent $1-\s$ could be reached with further elaboration.
In any case this is much weaker than the bound conjectured in [7] by the
author, namely that
for any given $\e > 0$ and fixed $\s$ satisfying $\hf < \s < 1$, one has
$$
\Z_2(\s + it) \;\ll_\e\; t^{{1\over2}-\s+\e}\qquad(t \ge t_0 > 0).
$$

\medskip
Both pointwise and mean square estimates for $\Z_2(s)$ may be used
to estimate $E_2(T)$. This connection is furnished by

\bigskip
THEOREM 2. {\it Suppose that for some $\rho \ge 0$ and $r \ge 0$ we have}
$$
\Z_2(\s + it) \;\ll_\e\; t^{\rho+\e},\quad
\int_1^T|\Z_2(\s + it)|^2\d t \;\ll_\e\; T^{1+2r+\e}\quad(\hf < \s \le 1).
\leqno(1.8)
$$
{\it Then we have}
$$
E_2(T) \;\ll_\e\; T^{{2\rho+1\over2\rho+2}+\e}, \quad
E_2(T) \;\ll_\e\; T^{{2r+1\over2r+2}+\e},\leqno(1.9)
$$
{\it and}
$$
\int_0^T|\zt|^8\d t \;\ll_\e\;T^{{4r+1\over 2r+1}+\e}.\leqno(1.10)
$$
 
\medskip\no
Note that from (1.2) with $\s = \hf + \e$ one can take in (1.8) $r = \hf$,
hence (1.9) gives
$$
E_2(T) \;\ll_\e\; T^{{2\over3}+\e}.\leqno(1.11)
$$
The bound (1.11) is, up to ``$\e$", currently the best known one
 (see [10] and [15], where $E_2(T) \ll T^{2/3}\log^8T$ is proved). Thus
any improvement of the existing mean square bound for $\Z_2(s)$ at
$\s = \hf + \e$ would result in (1.11) with the exponent strictly less than
2/3, which would be important.  Of course, if the first bound in (1.8) 
holds with some $\rho$, then trivially the second bound will hold with 
$r = \rho$. Observe that the known value $r = \hf$ and (1.10) yield 
$$
\int_0^T|\zt|^8\d t \;\ll_\e\; T^{3/2+\e},
$$
which is, up to ``$\e$", currently the best known bound for the eighth
moment (see [1, Chapter 8]), and ant value $r < \hf$ would reduce the
exponent 3/2 in the above bound.  

\bigskip
{\hh 2. The necessary lemmas}
\bigskip
\no
This section contains the lemmas needed for the proof of Theorem 1.
Let, as usual,
$\a_j = |\rho_j(1)|^2(\cosh\pi\k_j)^{-1}$, where
$\rho_j(1)$ is the first Fourier coefficient of  the Maass wave form
corresponding to the eigenvalue $\lambda_j$ to which the Hecke
$L$-function $H_j(s)$ is attached.

\medskip 
LEMMA 1. {\it We have}
$$
\sum_{K-G\le\k_j\le K+G} \a_j\H \;\ll_\e\; GK^{1+\e}
\quad(K^{\e}  \;\le\; G \;\le \; K).\leqno(2.1)
$$

\medskip\no
This result is proved by the author in [6]. Note that M. Jutila [12]
obtained
$$
\sum_{K-K^{1/3}\le\k_j\le K+K^{1/3}} \a_jH_j^4(\hf) \;\ll_\e\; K^{4/3+\e},
$$
but this result and (2.1) do not seem to apply each other. Both, however,
imply the hitherto sharpest bound for $H({1\over2})$, namely
$$
H_j(\hf) \;\ll_\e\; \k_j^{1/3+\e}.
$$
This bound is still quite far away from the conjectural bound
$$
H_j(\hf + it) \;\ll_\e\; (\k_j + |t|)^\e,
$$
which may be thought of as the analogue of the classical 
Lindel\"of hypothesis ($\zt \ll_\e |t|^\e$) for the Hecke series.

\medskip
LEMMA 2. {\it  Let $\xi\in (0,1)$ be a constant, and set}
$$ 
\psi(T)={1\over\sqrt{\pi}T^{\xi}}\int_{-\infty}^{\infty}
|\zeta(\hf + i(T+t))|^4 \exp(-(t/T^{\xi})^2)\d t.  \leqno(2.2)
$$
{\it Then we have}
$$ 
\psi(T)=I_{2,r}(T,T^{\xi})+I_{2,h}(T,T^{\xi})+I_{2,c}
(T,T^{\xi})+I_{2,d}(T,T^{\xi}).
\leqno(2.3)
$$  {\it
Here $I_{2,r}$ is an explicit main term, the contribution of
$I_{2,h}$ is small,
$$
I_{2,c}(T,T^{\xi })=\pi ^{-1}\int_{-\infty }^{\infty }
{|\zeta (\hf+ir)|^6
\over |\zeta (1+2ir)|^2}\Lambda (r;T,T^{\xi })\d r,
$$
$$ 
I_{2,d}(T,T^{\xi})=\sum_{j=1}^{\infty}\alpha _j  H_j^3(\hf)\Lambda
(\k_j;T,T^{\xi}),\leqno(2.4)
$$ 
where
$$
\eqalign{
\Lambda (r;T,T^{\xi}) &= \hf {\rm Re}\Biggl \{\left (1+
{i\over \sinh \pi r} \right )\Xi (ir;T,T^{\xi}) \cr&+ \left (1-
{i\over \sinh \pi r} \right )\Xi (-ir;T,T^{\xi}) \Biggr \}
\quad(r\in\RR)\cr }
\leqno(2.5)
$$ 
with
$$
\eqalign{
\Xi (ir;T,T^{\xi }) &={\G^2 ({1\over2}+ir)\over{\G (1+2ir)}}\int_0^{\infty }
(1+y)^{-{1\over2}+iT}y^{-{1\over2}+ir}\cr
& \times \exp\left(-\txt{1\over4} T^{2\xi}\log ^2(1+y)\right )
F(\hf +ir, \hf + ir;1+2ir;-y) \d y,\cr}
\leqno(2.6)
$$ 
and $F$ is the hypergeometric function}. 

\medskip\no
This fundamental result is
the spectral decomposition formula of Y. Motohashi (see [16, Section 5.1]
with $G = T^\xi$). The holomorphic part $I_{2,h}(T,T^\xi)$ is ``small",
namely by [16, Lemma 4.1] it is $\ll T^{-2}$. Motohashi (op. cit.)
gives an explicit evaluation of the main term $I_{2,r}(T,T^\xi)
\,(\ll \log^4T)$,
which can be used to show that the contribution from this function to
the relevant expression in Section 4 will be indeed absorbed by the
other terms. The structure of the continuous part $I_{2,c}(T,T^\xi)$
is similar in nature to (2.4), only the presence of integration instead
of summation over $\k_j$ makes this term less difficult to deal with than
(2.4).

\medskip
LEMMA 3. {\it The hypergeometric function, defined for $|z| < 1$ by}
$$
F(\alpha,\beta;\gamma;z) = 1+ \sum_{k=1}^\infty{(\alpha)_k(\beta)_k\over
(\gamma)_kk!}\,z^k\qquad ((\a)_k = \a(\a+1)\ldots\,(\a+k-1\,)),
$$
{\it satisfies}
$$
F(\a , \a ; 2\a ;z)=\left ({{1+\sqrt{1-z}}\over 2}\right )^{-2\a }
F\left (\a
, \hf ; \a +\hf ; \left ({1-\sqrt{1-z}}\over {1+\sqrt{1-z}}\right
)^2\right ).  \leqno(2.7)
$$

\medskip\no
This is a special case of the classical quadratic transformation formula
for the hypergeometric function (see e.g., [12, (9.6.12)]).

\medskip
LEMMA 4. {\it If $\,\Xi(ir;T,T^\xi)$ is defined by} (2.6), {{\it then}
$$
\Xi(ir;x,x^{\xi})\ll\cases{\phantom{|}1 &   $|r|\le x\log^2x$,\cr
|r|^{-C} & $|r| > x\log^2x$,\cr}\leqno(2.8)
$$ 
{\it which holds  uniformly for any fixed 
$\xi \in (0,1)$  and any fixed, large $C > 0$. We also have,
for $r\geq 2x^{1-\xi}\log ^5x$ and any fixed, large $D>0$, }
$$
\Xi (-ir;x,x^{\xi})\ll (rx)^{-D}.\leqno(2.9)
$$ 

\medskip\no 
The first bound in (2.8) (see [11, eq. (4.20)]) follows trivially from 
the defining relation (2.6) if the hypergeometric function is represented 
by the Gaussian integral formula (see e.g., [16, (6.17)-(6.18)]). The 
other bound in (2.8) and the one in (2.9) are contained in [16, Lemma 4.2].

\vfill\break
\bigskip
{\hh 3. The estimation of $\Z_2(s)$}
\bigskip
We are ready now to proceed with the estimation of $\Z_2(s)$. We
suppose that 
$\,\hf < \s_0 \le \s \le 1$, where $\s_0$ is fixed and let for some $C>1$
$$
T \le t \le 2T,\;s = \s + it,\;\hf < \s_0 \le \s \le 1,\;
  T^{1+\e} \le X \le T^C.\leqno(3.1)
$$
The reason for introducing $T$ is  for potential applications of our method
to mean square estimates of $\Z_2(s)$.
We start from the decomposition
$$
\eqalign{
&\Z_2(s) = \int_1^{2X}\rho(x)|\zx|^4x^{-s}\d x  \cr&
+ \int_X^{2Y}\s(x)(|\zx|^4 - \psi(x))x^{-s}\d x +
 \int_X^{2Y}\s(x)\psi(x)x^{-s}\d x
\cr& 
+ \int_Y^\infty \o(x)|\zx|^4x^{-s}\d x\cr&
= \Z_{12}(s) + \Z_{22}(s) + \Z_{32}(s) + \Z_{42}(s),\cr}\leqno(3.2)
$$
say. It is by introducing $\psi(T)$,  given by (2.2), that we are able
to exploit the spectral decomposition furnished by (2.3). 
We suppose that $0 < \xi \le \hf$, but eventually we shall take 
$\xi = \e$, namely arbitrarily small. This will follow, in the course of 
the estimation of $\Z_{32}(s)$, by an analysis similar to the one made 
in [7]. We suppose that $Y = Y(T,\s_0) \,(\ll T^C)$ is a large parameter
such that $Y>X$. 
The function $\rho(x) \,(\ge 0)$ is a smooth function supported in
$[1, 2X]$ such that $\rho(x) = 1$ for $1\le x \le X$, and
$\rho(x)$ monotonically decreases from 1 to 0 in $\,[X,\,2X]\,$.
The function $\s(x)$
is also a smooth non-negative function supported in $[X,\,2Y]\,$. 
We set $\s(x) = 1 - \rho(x)$ for
$X \le x \le Y$, and let $\s(x)$ monotonically decrease from 1 to 0 in
$[Y,\,2Y]\,$. Thus $\s(x)$ is supported in $\,[X,\,2Y]$, 
 $\s^{(\ell)}(X) = \s^{(\ell)}(2Y) = 0$ for 
$\ell = 0,1,2,\ldots\,$, and for $\ell\in\NN$ we have
$$
\s^{(\ell)}(x) = \cases{O_\ell(X^{-\ell})&$X\le x \le 2X$,\cr \cr 
0 &$2X \le x \le Y$,\cr \cr
O_\ell(Y^{-\ell})&$Y\le x \le 2Y$.\cr}\leqno(3.3)
$$
For $x \ge Y$ we set $\o(x) = 1-\s(x)$. Then we 
have  $\o^{(\ell)}(Y) = 0$ and $\o^{(\ell)}(x) \ll_\ell x^{-\ell}$ 
for $\ell = 0,1,2,\ldots\,$, and $\o'(x) = 0$ for $x \ge 2Y$.
This decomposition of $\Z_2(s)$ differs from the one that was made in [11].
Namely we have introduced here the parameters $X, Y$ and the smoothing 
functions $\rho, \s$ and $\o$.

Clearly the functions $\Z_{12}(s), \Z_{22}(s), \Z_{32}(s)$ are
entire functions for $s$ belonging to the region defined by (3.1).
The function $\Z_{42}(s)$ is initially defined for $\s > 1$, but we shall
presently see that it admits analytic continuation to the region
$\hf < \s_0 \le \s \le 1$,
and moreover its contribution (for $Y$ sufficiently large) 
will be negligible. To see this write (see (1.4))
$$
|\zx|^4 = Q_4(\log x) + E_2'(x),\quad Q_4(z) := P_4(z) + P_4'(z).\leqno(3.4)
$$
Then
$$
\Z_{42}(s) = \int_Y^\infty \o(x)(Q_4(\log x) + E_2'(x))x^{-s}\d x
\quad(\s > 1).\leqno(3.5)
$$
Integrating by parts we obtain
$$\eqalign{
\Z_{42}(s) &= {1\over s-1}\int_Y^\infty x^{1-s}
(\o(x)Q_4(\log x))'\d x\cr&
- \int_Y^\infty E_2(x)(\o'(x)x^{-s} - s\o(x)x^{-s-1})\d x
= \Z_{52}(s) + \Z_{62}(s),\cr}
$$
say. Since $\o'(x) \ll 1/x$, it follows from the mean square bound
for $E_2(T)$ (see e.g., [9]), by the Cauchy-Schwarz inequality 
for integrals, that $\Z_{62}(s)$ is regular for $\s > \hf$ and that
$$
\Z_{62}(s) \ll TY^{{1\over2}-\s}\log^CY
$$
holds for $s$ satisfying (3.1).   Now choose
$$
Y \;=\; T^{3\over2\s_0-1}.
$$
Then for $s$ satisfying (3.1) we have
$$
\Z_{62}(s) \ll TY^{{1\over2}-\s_0}\log^CY \ll T^{-{1\over2}}\log^{C}T,
$$
hence the contribution of $\Z_{62}(s)$ will be negligible. Repeated 
integration by parts gives, since $\o'(x)$ is supported in $[Y,\,2Y]$,
$$\eqalign{
\Z_{52}(s) = &\sum_{j=1}^5{1\over(s-1)^j}
\int_Y^{2Y}x^{1-s}\o'(x)Q_4^{(j-1)}(\log x)\d x\cr&
+ {1\over(s-1)^4}\int_Y^\infty x^{-s}\o(x)Q_4^{(4)}(\log x)\d x.\cr}
$$
Note that $Q_4^{(4)}(\log x) = C$, a constant, since $Q_4(z)$ is
a polynomial of degree four in $z$. Thus the last integral above becomes,
for $\ell \ge 2$,
$$\eqalign{
{C\over(s-1)^4}\int\limits_Y^\infty x^{-s}\o(x)\d x 
&= {C\over(s-1)^5(s-2)\ldots(s-\ell)}
\int\limits_Y^{2Y}x^{\ell-s}\o^{(\ell)}(x)\d x\cr&
\ll Y^{1-\s}T^{-4-\ell} \le Y^{1\over2}T^{-4-\ell} \ll T^{-{1\over2}}
\cr}
$$
on taking $\ell = \ell(\s_0)$ sufficiently large. The remaining integrals
with $Q_4^{(j-1)}(\log x)$ are treated in an analogous way. Integration by
parts is applied a large number of times, until each summand by trivial
estimation is estimated as $O(T^{-{1\over2}})$. Therefore the total
contribution of $\Z_{42}(s)$  will be negligible, as asserted.

Now we trivially estimate $\Z_{12}(s)$ by the fourth moment of $|\zt|$ as
$$
\Z_{12}(s) = \int_1^{2X}\rho(x)|\zx|^4x^{-s}\d x
\ll X^{1-\s}\log^4X,\leqno(3.6)
$$
for $s = \s + it$ and $\hf < \s_0 \le \s \le 1$. 

Next, the change of variable $t = \a x^\xi\log x$ 
in (2.2) gives
$$\eqalign{
\Z_{22}(s) = &-{1\over\sqrt{\pi}}
\int_{-\infty}^\infty \int_X^{2Y}
\left(|\z(\hf + i(x + \a x^\xi\log x))|^4 - |\zx|^4\right)\cr&\times
\s(x)\exp(-\a^2\log^2x)x^{-s}\log x\d x\d \a.\cr}
$$
The integral over $\a$
may be truncated with a negligible error at $|\a| = b$, with $b$
a small, positive constant. 
The relevant portion of $\Z_{22}(s)$ will be a multiple of
$$\eqalign{
\Z_{22}^*(s) :&= 
\int_{-b}^b\int_X^\infty
\left(E_2'(x + \a x^\xi\log x) - E_2'(x)\right)\cr&\times
\s(x)\exp(-\a^2\log^2x)x^{-s}\log x\d x\d \a, \cr} \leqno(3.7)
$$
where (3.4) is used. Namely, for $\xi = \e$,
the portion of $\Z_{22}(s)$ containing $Q_4$ 
makes a  total contribution that does not exceed the one in (3.6).

\medskip
In the $x$-integral in (3.7) we make the change of variable
$\tau = \tau (x,\alpha )= x+\alpha x^{\xi}\log x$. If $b $ is
sufficiently small,  then $\tau (x,\alpha )$ is
monotonically increasing as a function of $x$ for $x\geq 1$,  and
there is a monotonic inverse function $x=x(\tau ,\alpha )$. For
$x=x(\tau,\alpha )$, we have $x-\tau \ll(\log \tau )\tau ^{\xi}$, 
hence $x \asymp \tau$, and the implicit equation for $x$ shows that
$$ 
x(\tau ,\alpha )=\tau -  \alpha \tau ^{\xi}\log \tau  +O\left ( \tau
^{2\xi-1}\log ^2\tau
\right ).\leqno(3.8)
$$ 
Also we have
$$
{\partial x(\tau ,\alpha )\over\partial \tau } 
= 1-\alpha \tau^{\xi-1}(1+\xi\log
\tau)  + O(\tau ^{2\xi-2}\log ^2\tau).\leqno(3.9)
$$ 
For given positive $\alpha $, we combine the contributions of
$\alpha $ and $-\alpha $ in (3.7). In the respective integrals  we
put $\tau = \tau (x,\pm \alpha )$,  in the integral  involving
$E_2'(x)$ we simply change the notation $x$ to $\tau $, and we set
$$
G(u) \= \exp(-\a^2\log^2u)\cdot (1+\a u^{\xi-1}(1+\xi\log u))^{-1}
\s(u)\log u.
$$
Then in view of (3.9) it is seen that (3.7) becomes
$$\eqalign{
\Z_{22}^*(s) &= \int_0^b\int_{\tau(X,\a)}^\infty E_2'(\tau)
G(x(\tau,\a))(x(\tau,\a))^{-s}\d\tau\d\a\cr&
+ \int_0^b\int_{\tau(X,-\a)}^\infty E_2'(\tau)
G(x(\tau,-\a))(x(\tau,-\a))^{-s}\d\tau\d\a\cr&
- 2\int_0^b\int_X^\infty E_2'(\tau)\s(\tau)\tau^{-s}\exp(-\a^2\log^2\tau)
\log \tau\d\tau\d\a.\cr}\leqno(3.10)
$$
In the integrals with $G$ we replace the lower bounds of integration
in the $\tau$-integrals by $X$. Changing the variable back to
$x = x(\tau,\pm\a)$, using Lemma 1 and (1.14) it follows
that the total error made in this process will be $(\xi = \e)$
$$
\ll_\e\; X^{1-\s+\e}.\leqno(3.11)
$$
Now we write
$$
G(x(\tau,\pm\a)) = \exp(-\a^2\log^2\tau)\s(\tau)\log\tau + H(x(\tau,\pm\a)),
$$
say, where $H(x(\tau,\pm\a))$ is independent of $s$. 
Taking into account (3.8), it follows by using the mean value theorem that 
$$
H(x(\tau,\pm\a))   \ll \tau^{\xi-1}\log^3 \tau.
$$
The portion of the integrals in (3.10) with
$H(x(\tau,\pm\a))$, on changing the variable again to $x = x(\tau,\pm\a)$,
is estimated by Lemma 1. Its contribution does not exceed (3.11), hence
we are left wih
$$\eqalign{&
\Z_{22}^{**}(s) := \int_0^b\int_X^\infty E_2'(\tau)\Bigl\{
(x(\tau,\a))^{-s} \cr&+ (x(\tau,-\a))^{-s}
-2\tau^{-s}\Bigr\}\exp(-\a^2\log^2\tau)\s(\tau)\log\tau\d\tau\d\a.
\cr}\leqno(3.12)
$$
The expression in curly braces in (3.12) is expanded by
Taylor's formula at $\tau$. It becomes
$$\eqalign{&
(x(\tau,\a) + x(\tau,-\a) -2\tau)(-s\tau^{-s-1})
\cr&
\,+ {1\over2!}\left((x(\tau,\a) - \tau)^2 + (x(\tau,-\a) - \tau)^2\right)
{\partial\over\partial\tau}(-s\tau^{-s-1})
+ \ldots\,.
\cr}\leqno(3.13)
$$
Note that $|s|\tau^{-1} \;\ll\; TX^{-1} \;=\; T^{-\e}$
by (3.1), so that each time $\partial\over\partial\tau$
is taken in forming a new derivative in (3.13), its order will
decrease by a factor of $|s|\tau^{-1} \,(\ll\, T^{-\e})$.  We shall take 
sufficiently many terms in (3.13) in such a way that the error will 
make a negligible contribution (i.e., absorbed  by (3.11)), since (3.8) 
holds and $\xi = \e$. The expression in (3.13) will be
$$
\ll_\e T^2\tau^{\e-\s-2}.\leqno(3.14)
$$
Since (3.1) holds we obtain, by using the fourth moment for $|\zt|$,
$$
\Z_{22}(s) \;\ll_\e\; T^2  X^{\e-1-\s} \ll_\e X^{\e+1-\s}. \leqno(3.15)
$$

\medskip
We pass now to the contribution of $\Z_{32}(s)$. While the estimation
of $\Z_{12}(s)$ and $\Z_{22}(s)$ was essentially
elementary, it is the function $\Z_{32}(s)$ that is the most
delicate one in (3.2) and its treatment requires the
application of spectral theory, namely Lemma 2. 

\medskip
As discussed in [11] and in Section 2, after Lemma 2, the main contribution 
to $\Z_{32}(s)$ will come from $I_{2,d}$ in (2.4), namely from the 
discrete spectrum. Thus only this contribution  will be treated in
detail. It equals
$$
\int_X^{2Y} \s(x)I_{2,d}(x,x^{\xi})x^{-s}\d x  =\sum_{j=1}^{\infty}
\alpha_jH_j^3({\textstyle{1\over2}})
\int_X^{2Y}\s(x)\Lambda  (\k _j;x,x^{\xi})x^{-s} \d x,
\leqno(3.16)
$$ 
where (3.1) is assumed. The interchange of integration and summation
follows from (2.8) of Lemma 4, which ensures absolute convergence
on the right-hand side. Note that in place of the integral on the right-hand 
side of  (3.16)  we can consider
$$ 
X_r(s) \;:=\; \int_X^{2Y}\s(x)\Xi (-ir;x,x^{\xi})x^{-s}\d x
\qquad(r = \k_j),
\leqno(3.17)
$$
since the term  with $\Xi(ir;x,x^{\xi})$ (see (2.5)) has no 
saddle-point, and its estimation is less difficult. Next, note
that in view of (2.9) of Lemma 4 the sum in (3.16) can be
restricted to $\k_j \le T^{C_1}$ with a suitable constant
$C_1 = C_1(\s_0,\xi)$, since the tails of the series will make
a negligible contribution.

We make the  change of variable $y = z/x$ in the $\Xi$-integral (see (2.6))
in (3.17). This is done to regulate the 
location of the corresponding saddle point, similarly as in [7] and [11]. 
After the change of variable the integral  $X_r(s)$  becomes
$$ 
{\G^2 ({1\over2}-ir)\over{\G (1-2ir)}}
\int_{X}^{2Y} \s(x) x^{-{1\over2}+ir-s}L^*(r;x) \d x,\leqno(3.18)
$$ 
where 
$$
\eqalign{& L^*(r;x):= \int_0^{\infty}z^{-{1\over2}-ir} 
\left (1+{z\over x}\right )^{-{1\over2}+ix}
\cr &\times\exp \left
(-\txt{1\over4} x^{2\xi}\log ^2 \left (1+{z\over x}\right )\right )
F\left (\hf -ir, \hf  -ir; 1-2ir;-{z\over x}\right )
\d z.\cr}\leqno(3.19)
$$
In (3.19) we consider separately the ranges $ z/x \leq x^{-\delta}$ and 
$z/x > x^{-\delta}$ for a sufficiently small, fixed $\delta > 0$. 
In the latter range, the exponential factor  is 
$\ll x^{-A}$ for any fixed positive $A$ provided that $\xi > \delta$, 
which we may assume, and thus the total
contribution of the range $z/x > x^{-\delta}$ in (3.19) is negligible. 
Therefore so far we have reduced the problem to the estimation of 
a finite sum over $\k_j$ in (3.16) and a finite $z$-integral in (3.19).

\medskip
In the range $ z/x \leq x^{-\delta}$ in (3.19)
we transform the hypergeometric function  by (2.7) of Lemma 3, noting
that the new hypergeometric series converges rapidly. Namely in the
series expansion for the hypergeometric function (see Lemma 3)
we can take a finite number of terms so that the tails of
the series will make a negligible contribution. Each term, since in our
case $\alpha = \hf - ir, \beta = \hf, \gamma = 1 - ir$, will yield
similar expressions, and each contribution will be smaller than the
one coming from the preceding term.
Therefore the most significant term in the above series expansion
will be simply the leading term 1, so it suffices to
consider its contribution. Then the
essential part of $L^*(r;x)$, say $L(r;x)$, takes the form
$$
\eqalign{& L(r;x) :=\, 2^{1-2ir}\int_0^{x^{1-\delta}}
z^{-{1\over2}} \left (1+{z\over x}\right )^{-{1\over2}}
\cr &\times\exp\left(
-{\txt{1\over4}} x^{2\xi}\log ^2 \left (1+{z\over x}\right )\right )
\left (1+\sqrt{1+ {z\over x}}\,\right )^{-1}
e^{i\f(z)} \d z,\cr}\leqno(3.20)
$$
where
$$
\f(z) = \f(r,x;z) := -r\log z + x\log\left (1+{z\over x}\right )
+2r\log\left(1 + \sqrt{1+{z\over x}}\,\right).
\leqno(3.21)
$$
The integral in (3.20) can  be approximately evaluated by the saddle 
point method (see e.g., [2, Chapter 2]). The main contribution to
the integral  comes from the  saddle point $z_0$ satisfying
$\f'(z_0) = 0$. We have
$$
\f'(z) = {\partial\f\over \partial z} = -{r\over z} + {x\over x+z} +
{r\over x\left(\sqrt{1+{z\over x}} + 1 + {z\over x}\right)},\leqno(3.22)
$$ 
and
$$
\f''(z) = {r\over z^2} - {x\over(x+z)^2}
- {r\over x^{2}}\cdot{{1+{1\over2\sqrt{1+{z\over x}}}\over
\left(\sqrt{1+{z\over x}} + 1 + {z\over x}\right)^2}}.\leqno(3.23)
$$
This gives
$$
z_0 = r\left(1 + {r\over2x} + {r^2\over8x^2} 
+ O\left({r^3\over x^3}\right)\right),\leqno(3.24)
$$
and the error term in (3.24) admits an asymptotic expansion in powers of 
$r/x$, since the relevant range is $r/x \ll T^{-\e}$, in view of
(3.1) and (3.32). Similarly we find that
$$
{\left({z_0\over x}\right)}'
= {\partial\over \partial x}\left({z_0\over x}\right)
= -{r\over x^2} - {r^2\over x^3} + O\left({r^3\over x^4}\right),\leqno(3.25)
$$
$$
z_0' = r\left( - {r\over2x^2} - {2r^2\over8x^3} 
+ O\left({r^3\over x^4}\right)\right)
,\;\f''(z_0) = {1\over r}\left(1+ O\left({r\over x}\right)\right),
\leqno(3.26)
$$
where again the error terms admit an asymptotic expansion in powers of $r/x$.
A calculation then shows that we obtain
$$
{\partial \f(z_0)\over\partial x} =  - {r^3\over24x^3} +
O\left({r^4\over x^4}\right).\leqno(3.27)
$$

As already asserted the main contribution to the integral in (3.20) 
comes from the  saddle point $z_0$, and equals a multiple of
$$
C(z_0)(\f''(z_0))^{-{1\over2}}e^{i\f(z_0)},\leqno(3.28)
$$
where
$$\eqalign{
C(z) \;&:=\; C(\xi,x;z) = z^{-{1\over2}} 
\left (1+{z\over x}\right )^{-{1\over2}}\cr &\times\exp \left
(-\txt{1\over4} x^{2\xi}\log ^2 \left (1+{z\over x}\right )\right )
\left(1+ \sqrt{1+{z\over x}}\,\right)^{-1}.\cr}\leqno(3.29)
$$
From (3.24) and (3.26) we have $(z_0\f_0''(z_0))^{-1/2} \sim 1$. Hence
inserting (3.28)-(3.29) in (3.18) it is seen that the main contribution
will be (by using Stirling's formula to simplify the gamma-factors)
a multiple of
$$
\eqalign{&
r^{-1/2}\int_X^{2Y}\s(x)x^{-\hf + ir-s}\left(z_0\f''(z_0)\right)^{-1/2}
\left(1 + {z_0\over x}\right)^{-1/2}\times\cr&
\times \exp\left(-{\txt{1\over4}}x^{2\xi}\log^2
\left(1 + {z_0\over x}\right)\right)
\left(1+\sqrt{1 + {z_0\over x}}\,\right)^{-1}e^{i\f(z_0)}\d x \cr&
= r^{-1/2}\int\limits_X^{2Y}\s(x)x^{-\hf  + ir-s}g(x)
\exp\left(-ir\log r+ {ir^3\over 48x^2}+ih(r,x)\right)\d x,\cr}\leqno(3.30)
$$
say, where in view of (3.27) we have
$$
{\partial^\ell h(r,x)\over\partial x^\ell}
\asymp_\ell {r^4\over x^{3+\ell}},\quad {\d^\ell g(x)\over\d x^\ell}
\ll_\ell x^{-\ell}\quad(\ell = 0,1,2,\ldots\,).
$$
From the term
$
\exp\left(-{\txt{1\over4}}x^{2\xi}\log^2
\left(1 + {z_0/x}\right)\right)
$
it transpires that the range $r \ge x^{1-\xi}\log x$ will make a
negligible contribution. Similarly if
$$
|r - t| > r^\a > {r^3t^\e\over x^2},\leqno(3.31)
$$
repeated integration by parts shows that the contribution is negligible.
But as $r \le x^{1-\xi}\log x$, (3.31) will hold for $(3-\a)(1-\xi)<2$,
namely for
$$
\a > 3 - {2\over 1-\xi}.
$$
If $0 < \xi < 1/3$, as we assume, then it follows that $0 < \a < 1$,
hence only the range $|r-t| < r^\a$ is relevant, which gives
$$
C_1T \le r = \k_j \le C_2T\qquad(0 < C_1 < C_2)\leqno(3.32)
$$
for suitable constants $C_1,C_2$.

Repeated integration by parts in (3.30) shows then that, after each 
integration, the order of the integrand is lowered by a factor of
$$
{x\over1+|t-r|}\left({1\over x} + {r^3\over x^3}\right).\leqno(3.33)
$$
It follows that, if $x \ge r^{3/2} \asymp T^{3/2}$, then the above
expression is $\ll T^{-\e}$ for $|r-t| > T^\e$, and hence its total
contribution is negligible.  If $|r-t| \le T^\e$, then by trivial
estimation and Lemma 1 the total contribution will be
$$
\ll t^{-1/2}\sum_{|r-t| \le T^\e}\a_j\H X^{{1\over2}-\s} \ll_\e
T^{{1\over2}+\e}X^{{1\over2}-\s} \ll_\e t^{1-\s+\e},
$$
hence absorbed by the right-hand side of (1.7). In view of (3.33)
this shows that the relevant range for $x$ and $r$ is
$$
T^{1+\e} \le X \le x \le T^{3/2},\; T^\e < |r-t| < T^{3+\e}X^{-2}.
\leqno(3.34)
$$

\medskip
The main contribution in (3.30) will, as in the previous case, come also
from the saddle point. If we write the integral as 
$$
r^{-1/2}\int_X^{2Y}\s(x)g(x)x^{-{1\over2}-\s}e^{iH(r,x)}\d x,
$$
where
$$
H(r,x)  :=  (r-t)\log x + {r^3\over 48x^2} + h(r,x),
$$
then in this case the saddle point $x_0$ is
the solution (in $x$) of $H'(r,x) = 0$, so that
$$
x_0 = \sqrt{{r^3\over24(r-t)}}\left(1 + O\left({r\over x_0}\right)\right),
\leqno(3.35)
$$
where the $O$-term admits an asymptotic expansion. We have also
$$
H''(r,x_0) \sim 48(r-t)^2r^{-3},
$$
and we suppose $G < r - t \le 2G,$ $G = 2^kT^\e,\,k\in\NN$. Then  it is seen 
that the contribution coming from the ranges
$$
x \ge C_1T^{3/2}G^{-1/2}, \quad x \le C_2T^{3/2}G^{-1/2}
$$
for sufficiently large $C_1$ and sufficiently small $C_2 > 0$ will be
negligible. There remains a multiple of ($S(x) := \s(x)g(x)x^{-1/2-\s}$)
$$
\eqalign{
\sum &:= \sum_{t+G< \k_j\le t+2G}\a_j\H \k_j^{-1/2}S(x_0)
(H''(\k_j,x_0))^{-1/2}e^{iH(\k_j,x_0)}\cr&
\ll \sum_{t+G< \k_j\le t+2G}\a_j\H\k_j^{-1/2}x_0^{-1/2-\s}x_0^2
\k_j^{-3/2}\cr&
\ll T^{-2}\sum_{t+G< \k_j\le t+2G}\a_j\H(T^{3/2}G^{-1/2})^{3/2-\s}\cr&
\ll_\e T^{\e-1}GT^{9/4-3\s/2}G^{\s/2-3/4} = T^{5/4-3\s/2+\e}G^{1/4+\s/2},
\cr}
$$
where Lemma 1 was used.
From (3.34) it follows that $G \ll T^{3+\e}X^{-2}$, hence we obtain
$$
\sum \;\ll_\e T^{2+\e}X^{-1/2-\s},
$$
and the same bound will hold for $\Z_{32}(s)$. Hence finally, in view 
of (3.6) and (3.15), we obtain
$$
\Z_2(\s + it) \ll_\e t^\e(X^{1-\s} + t^2X^{-1/2-\s}) \ll_\e 
t^{{4\over3}(1-\s)+\e}\quad(\hf < \s \le 1)
$$
for $X = t^{4/3}$. This completes the proof of Theorem 1 .

\bigskip
{\hh 4. Proof of Theorem 2}

\bigskip\no
The proof follows the analysis outlined in [7], where the first
bound in (1.9) was mentioned.
From the defining relation (1.4) it is not difficult to obtain
(see e.g., [4, (5.3)]) that $(C_1, C_2 > 0,\;1 \ll H \le {\txt{1\over4}}T)$,
$$
E_2(T) \le C_1H^{-1}\int_T^{T+H}E_2(x)f(x)\d x + C_2H\log^4T,\leqno(4.1)
$$
where $f(x)\;( > 0)\,$ is a smooth function supported in $\,[T,\,T+H]\,$,
such that $f(x) = 1$ for $T + {1\over4}H \le x \le T + {3\over4}H$.
Taking account that $\Z_2(s)$ is the (modified) Mellin transform of
$|\zx|^4$, it follows by the Mellin inversion formula that (see [7, (3.27)])
$$
|\zx|^4 \= {1\over2\pi i}\int_{\cal L}{\cal Z}_2(s)x^{s-1}\d s
+ Q_4(\log x) \qquad(x > 1), \leqno(4.2)
$$
where we have set, as in (3.4), 
$Q_4(\log x) \= P_4(\log x) + P'_4(\log x)$ and $\cal L$ denotes the
line $\Re s = 1 + \e$ with a small indentation to the left of $s=1$. 
If we integrate (4.2) from $x= 1$ to $x = T$ and
take into account the defining relation  of $E_2(T)$, we shall obtain
$$
E_2(T) \= {1\over2\pi i}\int_{\cal L}{\cal Z}_2(s){T^s\over s}\d s
+ O(1)\qquad(T > 1).\leqno(4.3)
$$
Then from (4.1) and (4.3) we have by Cauchy's theorem ($\hf < c < 1,\,T>1$)
$$
E_2(T) \le {C_1\over2\pi iH}\int_{(c)}{{\cal Z}_2(s)\over s}
\int_T^{T+H} f(x)x^s\d x\,\d s + C_2H\log^4T,\leqno(4.4)
$$
and we also have an analogous lower bound for $E_2(T)$.
Since $f^{(r)}(x) \ll_r H^{-r}$ it follows that the $s$--integral
in (4.4) can be truncated at $|\Im{\rm m}\, s|= T^{1+\e}H^{-1}$ with
a negligible error, for any $c$ satisfying $\hf < c < 1$.
We take $c = \hf + \e$ and use the first bound in (1.8) to obtain
$$\eqalign{
E_2(T) &\ll_\e \int_1^{T^{1+\e}H^{-1}}t^{\rho-1+\e}T^{{1\over2}+\e}\d t 
+ HT^\e \cr&\ll_\e T^\e(T^{{1\over2}+\rho}H^{-\rho} + H) 
\ll_\e T^{{2\rho+1\over2\rho+2}+\e}\cr}
$$
with the choice $H = T^{2\rho+1\over2\rho+2}$. This proves the first
part of Theorem 2. To prove the second, we proceed similarly, but use
the Cauchy-Schwarz inequality and the second bound in (1.8). We have
$$\eqalign{
E_2(T) &\ll_\e \int_1^{T^{1+\e}H^{-1}}|\Z_2(\hf+\e+it)|t^{-1}
T^{{1\over2}+\e}\d t + HT^\e \cr&\ll_\e 
T^{{1\over2}+\e}
\left(\int_1^{T^{1+\e}H^{-1}}|\Z_2(\hf+\e+it)|^2t^{-1}\d t\right)^{1/2} 
+ HT^\e\cr&
\ll_\e T^\e(T^{{1\over2}+r}H^{-r} + H) \ll_\e T^{{2r+1\over2r+2}+\e}
\cr}
$$
with $H = T^{2r+1\over2r+2}$. 

Finally to prove (1.10), note that by [8, eq. (4.9)] we have
($c_2(\hf+\e) = 1 + 2r$ in this notation)
$$
\int_0^T|\zt|^{4+4A}\d t \ll_\e T^{A+\e}\qquad(1 \le A \le 2).\leqno(4.5)
$$
Let $0 \le C \le 8$ be a constant. Then, for $p > 0, q > 0, 1/p+1/q=1$,
H\"older's inequality for integrals gives
$$
\eqalign{&
\int_0^T|\zt|^8\d t = \int_0^T|\zt|^C|\zt|^{8-C}\d t\cr&
\le \left(\int_0^T|\zt|^{Cp}\d t\right)^{1/p}
\left(\int_0^T|\zt|^{(8-C)q}\d t\right)^{1/q}.\cr}
$$
Now choose $p,q$ and $C$ so that
$$
Cp = 4,\;(8-C)q = 4 + 4A,\;{1\over p} + {1\over q} = 1.
$$
Using the fourth moment and (4.5) the above inequality gives then
$$
\int_0^T|\zt|^8\d t \ll_\e T^{{2A-1\over A}+\e} = T^{{4r+1\over 2r+1} +\e},
$$
since $A  = 1 + 2r$. This completes the proof of Theorem 2.

In concluding let us remark that $r = 0$ in (1.8) was conjectured by
the author in [7]. This is a very strong conjecture since it gives,
by (1.6) and (1.9),
$$
\int_0^T|\zt|^8\d t \ll_\e T^{1+\e},\quad E_2(T) \ll_\e T^{{1\over2}+\e},
\leqno(4.5)
$$
and both of these bounds are, up to``$\e$", best possible. Namely one has
(see e.g., [1] and [16])
$$
\int_0^T|\zt|^8\d t \gg T\log^{16}T,\quad E_2(T) = \Omega_\pm(\sqrt{T}).
$$
So far it is not known whether either of the bounds in (4.5) implies the
other one. However, it is known that either of them
implies  (see [2, Theorem 1.2 and Lemma 4.1])  
the hitherto unproved bound
$$
\zt \;\ll_\e\;|t|^{{1\over8}+\e}.
$$


\bigskip
\centerline{\hh REFERENCES}
\bigskip

\item {[1]} A. Ivi\'c,  The Riemann zeta-function, {\it John Wiley and
Sons}, New York, 1985.

\item {[2]} A. Ivi\'c,  Mean values of the Riemann zeta-function, LN's
{\bf 82}, {\it Tata Institute of Fundamental Research}, Bombay,  1991
(distr. by Springer Verlag, Berlin etc.).

\item{[3]}  A. Ivi\'c,  On the fourth moment of the Riemann
zeta-function, {\it Publs. Inst. Math. (Belgrade)} {\bf 57(71)}
(1995), 101-110.

\item{[4]} A. Ivi\'c,  The Mellin transform and the Riemann
zeta-function,  {\it Proceedings of the Conference on Elementary and
Analytic Number Theory  (Vienna, July 18-20, 1996)},  Universit\"at
Wien \& Universit\"at f\"ur Bodenkultur, Eds. W.G. Nowak and J.
Schoi{\ss}engeier, Vienna 1996, 112-127.

\item{[5]} A. Ivi\'c, On the error term for the fourth moment of the
Riemann zeta-function, {\it J. London Math. Soc.} {\bf60}(2)(1999), 21-32.

\item{[6]} A. Ivi\'c, On sums of Hecke series in short intervals,
{\it J. de Th\'eorie des Nombres Bordeaux} {\bf13}(2001), 1-16.

\item{[7]} A. Ivi\'c, On some conjectures and results for the
Riemann zeta-function, {\it Acta Arith.} {\bf109}(2001), 115-145. 

\item{[8
]} A. Ivi\'c, Some mean value results for the Riemann
in: Jutila Matti (ed.) et. al. Number Theory. Proc.
Turku Symposium 1999, Berlin, de Gruyter 2001, 145-161.
zeta-function,

\item{ [9]} A. Ivi\'c and Y. Motohashi,  The mean square of the error
term for the fourth moment of the zeta-function,  {\it Proc. London
Math. Soc.} (3){\bf 66}(1994), 309-329.

\item {[10]} A. Ivi\'c and Y. Motohashi,  The fourth moment of the
Riemann zeta-function, {\it J. Number Theory} {\bf 51}(1995), 16-45.

\item {[11]} A. Ivi\'c, M. Jutila and Y. Motohashi, The Mellin
transform of powers of  the Riemann zeta-function,  {\it Acta Arith.}
{\bf95}(2000), 305-342.

\item{[12]} M. Jutila,  The fourth moment of central values of Hecke series,
in: Jutila Matti (ed.) et. al. Number Theory. Proc.
Turku Symposium 1999, Berlin, de Gruyter 2001, 167-177.

\item{[13]} N.N. Lebedev, Special functions and their applications,
{\it Dover}, New York, 1972.

\item{ [14]} Y. Motohashi,   An explicit formula for the fourth power
mean of the Riemann zeta-function, {\it Acta Math. }{\bf 170}(1993),
181-220.

\item {[15]} Y. Motohashi,  A relation  between the Riemann
zeta-function and the hyperbolic Laplacian, {\it Annali Scuola Norm.
Sup. Pisa, Cl. Sci. IV ser.} {\bf 22}(1995), 299-313.

\item {[16]} Y. Motohashi,  Spectral theory of the Riemann
zeta-function, {\it Cambridge University Press}, Cambridge, 1997.



\bigskip
\parindent=0pt

\cc
Aleksandar Ivi\'c \par
Katedra Matematike RGF-a
Universiteta u Beogradu\par
Dju\v sina 7, 11000 Beograd,
Serbia (Yugoslavia)\par
{\sevenbf e-mail: aivic@matf.bg.ac.yu, 
aivic@rgf.bg.ac.yu}

\bye